# ASPECTS ANALYTIQUES DANS LA MATHEMATIQUE DE SHARAF AL-DÎN AL-TÛSÎ

Nicolas FARÈS[*]

## 0. INTRODUCTION.

Dans le présent chapitre, nous nous proposons d'aborder le caractère analytique de la mathématique du "Traité des équations" de Sharaf al-Dîn al-Tûsî. Nous insistons surtout sur l'impossibilité de classer dans le seul domaine de la géométrie euclidienne certaines techniques fondamentalement algébriques de cette œuvre.

Le "*Traité*" est connu grâce à R. Rashed qui, après en avoir exposé en 1974 les principaux résultats [1], l'a établi, traduit et commenté en 1986, dans un livre de deux tomes [2], avec l'ensemble des œuvres d'al-Tûsî connues de nos jours. Rashed, tout en insistant sur l'unité du "Traité", l'a divisé en deux parties, incluses chacune dans un de ces deux tomes. Bien qu'elle puisse être motivée par des soucis d'ordre pratique inhérents à l'édition, cette subdivision correspond bien à une différenciation d'ordre mathématique (1).

La publication de cette œuvre "renouvelle profondément la connaissance que nous avions des mathématiques arabes consacrées aux équations algébriques" [4]. A ce titre, elle attire l'attention des historiens des mathématiques. Il n'est pas habituel, en effet, de rencontrer, réunies dans un même traité, des disciplines variées qui sont: l'algèbre, la géométrie, "l'analyse" et le calcul numérique; dans cette œuvre d'al-Tûsî, celles-ci sont là pour l'élaboration d'une théorie des équations du troisième degré; laquelle fut le sujet d'intérêt de générations de mathématiciens tout au long de la période séparant Al-Khawârizmî d'Al-Khayyâm. La concertation de ces disciplines mathématiques autour d'un thème unique, les nouveautés introduites par al-Tûsî tant au niveau conceptuel qu'à celui des méthodes et des techniques calculatoires, font du "Traité", selon Rashed, "la plus importante œuvre algébrique écrite en langue arabe " [2. t. 1. p. V; préface arabe, li.

[*] Nicolas FARÈS
Département des mathématiques – Faculté des Sciences – Université de Reims.
Département des mathématiques – Faculté des Sciences – Université Libanaise.
Equipe d'Etude et de Recherche sur la Tradition Scientifique Arabe, CNRS-Liban.

5]; déclaration de taille, vu la place qu'occupait cette discipline dans la littérature mathématique arabe tout au long de l'époque médiévale. Fruit d'une synthèse de deux courants algébriques: l'algèbre "arithmétique" et l'algèbre "géométrique", la mathématique du "*Traité*" est, selon C. Houzel, une "mathématique de très haut niveau pour l'époque " [3. p. 62].

Outre la variété et la richesse du "Traité", la publication commentée de cette œuvre ne manque pas de suciter l'intérêt du chercheur, pour des raisons parmi lesquelles nous citons:

1) Le style d'al-Tûsî ne laisse apparaître explicitement que sa synthèse et presque nullement son analyse. Il montre, généralement, qu'une condition est suffisante pour aboutir à un résultat, sans toutefois indiquer le processus analytique qui l'a conduit à la découverte de cette condition. Tel est, par exemple, le cas quand il avance que certaines expressions polynomiales f(x) passent par un maximum quand x est racine d'une équation équivalente à f'(x) = 0. Un tel style pourrait offrir des terrains de recherche en vue d'une "reconstitution" de l'analyse manquante.

2) Dans les premiers chapitres, consacrés à l'analyse de l'œuvre d'al-Tûsî et à son commentaire, R. Rashed donne une étude approfondie du "*Traité*" sur les deux plans historique et mathématique. Les commentaires accompagnant la traduction constituent un outil indispensable à la compréhension du texte, souvent extrêmement difficile (ici nous ne faisons pas seulement allusion à la reconstitution des tableaux revenant aux algorithmes numériques complètement supprimés par le copiste, ni au langage du texte). "L'Introduction" de Rashed met l'accent sur le style synthétique d'al-Tûsî; elle propose certaines interprétations et conjectures, invitant explicitement à la poursuite des recherches historiques afin de réaliser une avancée dans l'explication de certains points de cette œuvre mathématique.

Notre démarche s'inscrit dans le cadre suggéré par cette "Introduction": elle consiste à étudier de près la "deuxième partie" du "*Traité*" en vue de chercher, dans les détails les plus subtils de cette mathématique, des éléments qui puissent s'ajouter aux données qui la situent et la caractérisent. Dans cette optique, nous jugeons bon de nous arrêter sur des points bien précis:

1/ la détermination algébrique des racines des équations de la deuxième partie du "*Traité*"[1];

2/ le rôle du "maximum" (concept et valeur effective) dans la démonstration de l'existence, la détermination et le calcul des racines des équations de cette partie;

3/ le rôle secondaire de la géométrie euclidienne dans cette partie.

**1. METHODE DE DETERMINATION DE LA PETITE RACINE.**

Le style synthétique du "Traité" implique, pour interpéter certains résultats, le recours à des "reconstitutions" analytiques de certains détails. De telles démarches, en l'absence de preuves sûres, pourraient engendrer des constatations susceptibles de varier d'une étude à une autre. Citons, à titre d'exemple, la reconstitution de la méthode relative à la détermination de la "petite racine" de l'équation (21) (pour la numérotation des équations, cf. note(2)); l'importance de cette détermination (cf. 1.1. f, plus bas) exige qu'elle soit examinée attentivement. Pour cela et, par souci de clarté, rappelons le schéma de la démarche d'al-Tûsî dans la résolution de cette équation et des quatre autres (22),..,(25) qui la suivent. Ces cinq équations cubiques, classées et traitées en dernier par al-Tûsî, sont celles qui peuvent ne pas avoir de racines (réelles positives). Elles occupent la deuxième partie du "Traité".

**1. 1. Schéma de la méthode de résolution des équations (21)..(25).**

Il s'agit des équations qui peuvent ne pas avoir de racine positive; elles sont les suivantes:

**(21)** $x^3 + c = ax^2$**,**
**(22)** $x^3 + c = bx$**,**
**(23)** $x^3 + ax^2 + c = bx$**,**
**(24)** $x^3 + bx + c = ax^2$**,**
**(25)** $x^3 + c = ax^2 + bx$**.**

Remarquons que dans chacune d'elles, $c$ et $x^3$ sont dans le même membre[2]. Le schéma de leur résolution par al-Tûsî est uniforme (dans ses grandes lignes) et se trouve

---

[1] Voir la note suivante.
[2] La numérotation des équations est celle adoptée par R. Rashed conformément à l'ordre de leur traitement dans le *Traité* d'al-Tûsî. Les équations numérotées de 1 à 20 occupent la première partie du Traité; chacune d'elles (sauf la 9e qui est une équation du 2e degré) admet toujours une racine réelle positive. Les équations

exposé dans [2. t. 1. pp. XXII - XXIII]; nous pouvons le décrire, en détail, de la façon suivante:

1) Al-Tûsî considère chacune de ces équations et la met sous la forme

$$(E) \qquad f(x) = c;$$

par exemple, l'équation (23) est mise sous la forme $bx - ax^2 - x^3 = c$; (cette écriture n'est pas conforme à la tradition algébrique: *al-jabr* consiste à ne laisser aucun membre de l'équation contenir des termes soustractifs).

2) Il détermine des conditions sur $x$, pour que l'équation $f(x) = c$ soit "possible"; ce qui revient à considérer le domaine **D** de définition de $f$ (l'ensemble des réels positifs $x$ tels que $f(x) > 0$).[3]

3) Il considère alors une racine $x_0$ d'une équation du $2^e$ degré équivalente à l'équation $f'(x) = 0$[4] et il annonce que $f(x)$ prend une valeur maximum[5] $f(x_0)$ pour $x = x_0$. Il prouve ensuite que $f(x_0)$ est le maximum de $f(x)$ en démontrant que:

$$\begin{cases} x < x_0 \Rightarrow f(x) < f(x_0) \\ x > x_0 \Rightarrow f(x) < f(x_0) \end{cases}$$

4) Il affirme que :

    4.1: si $c > f(x_0)$, l'équation n'a pas de racines (réelles positives);

    4.2: si $c = f(x_0)$, l'équation admet une seule racine (double);

    4.3: si $c < f(x_0)$, elle admet deux racines $x_1$ et $x_2$ telles que:

$$0 < x_1 < x_0 < x_2.$$

---

21-25 "*pouvant être impossibles*", occupent la deuxième partie. R. Rashed édite, traduit et commente le *Traité* d'al-Tûsî, dans un livre de deux tomes dont chacun correspond à une partie du *Traité*.

[3] **D** est bien déterminé pour les équations (21) et (22) où on a respectivement **D** = ]0, $a$[ et **D** = ]0, $\sqrt{b}$ [; pour les équations (23) et (24), al-Tûsî ne donne que des conditions nécessaires pour que $f(x)$ soit > 0. De telles conditions sont absentes pour l'équation (25). En revanche, pour cette dernière, ainsi que pour l'équation (24), il présente une étude qui revient à localiser les deux racines positives.

[4] Par exemple, dans le cas de l'équation (23), il s'agit de l'équation: $x^2 + 2\dfrac{a}{3}x = \dfrac{b}{3}$, i.e. $\dfrac{1}{3}f'(x) = 0$. On peut interpréter la présence du facteur multiplicatif $\dfrac{1}{3}$ par l'usage courant de l'époque qui consiste à présenter une équation sous sa forme normale (ou canonique, avec 1 comme coefficient du terme du plus haut degré en x).

[5] Il s'agit plus précisément, du "nombre maximum" ou *plus grand nombre*"; "*le nombre*" pour al-Tûsî étant la constante $c$ dans l'équation $f(x) = c$. Si la valeur de $c$ dépasse une certaine valeur $f(x_0)$ alors le problème est impossible; $x_0$ étant une racine de $f'(x) = 0$.

A la suite de cette affirmation et sous l'hypothèse $c < f(x_0)$, il effectue les étapes 5), 6), 7) et 8) suivantes:

5) Il détermine $x_2$ par un changement de variable affine: $x = x_0 + X$, qui ramène l'équation (E) à une équation du type (15) en $X$:

$$(15) \quad X^3 + aX^2 = c,$$

qui a une seule racine -réelle positive- et qui a été résolue dans la première partie de son "*Traité*".

6) Il détermine $x_1$ (sauf pour l'équation (21)) par un changement de variable affine: $x = x_0 - X$, qui ramène (E) à une équation du type (21) en $X$: ($aX^2 - X^3 = c$). Si $X_1$ est la petite racine de cette nouvelle équation, alors $x_0 - X_1$ est celle de (E). R. Rashed remarque [2. t. 2. p. XXIV] qu'un tel changement de variable appliqué à l'équation (21), la ramène à une équation du même type, d'où le traitement exceptionnel de cette équation par al-Tûsî, qui ramène sa solution à celle d'une équation du type (7):

$$(7) \quad x_2 + ax = c$$

déjà résolue dans la première partie [2, t. 1, p. CXLV].

7) Il calcule $x_2$ numériquement, moyennant l'algorithme propre à l'équation (15) (déjà résolue).

8) Il calcule $x_1$ numériquement, moyennant un algorithme propre à l'équation (21).

Notons toutefois que les étapes 1) et 2) forment pratiquement une même étape qu'on pourra noter [1), 2)] et que, seulement dans l'équation (22), l'étape 6) précède l'étape 5).

**1. 2. Détermination de $x_1$.**

Considérons maintenant la détermination de la petite racine $x_1$ de l'équation (21) ( $x_2$ étant déjà déterminée, - voir 5) et 7), plus haut -). On est dans le cas:

$$0 < x_1 < x_0 = 2a/3 < x_2 < a.$$

Al-Tûsî considère alors l'équation du type (7):

$$X^2 + (a - x_2).X = x_2(a - x_2)$$

admettant une seule racine positive *X* et il démontre que $x_1 = (a - x_2) + X$, sans indiquer les raisons de cette démarche purement synthétique (la partie analytique de son raisonnement étant absente comme c'est souvent le cas).

La reconstitution de J. P. Hogendijk dans [5. p. 78. li. 7 .. 26], d'ailleurs vraisemblable, ne nous paraît pas conduire nécessairement à une interprétation géométrique du raisonnement d'al-Tûsî. Reprenons ses notations qui, ici, sont aux nominations près, celles du "*Traité*", et considérons la figure qu'il utilise, tirée de [5. p. 76. fig. 2]:

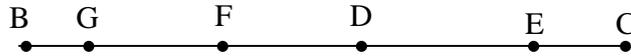

où $BC = a$, $BD = x_0 = 2a/3$, $BE = x_2$, $BF = x_1$, $BG = EC = a - x_2$. (La petite racine de l'équation (21), que nous notons $x_1$ est notée $x_2$ dans [5], tandis que la grande est notée $x_1$). Avec nos notations qui sont celles de [2], on a:

$$c = ax_2^2 - x_2^3 = BE^2.CE = ax_1^2 - x_1^3 = BF^2.CF \ ;$$

ce qui donne, en retranchant $BE^2.CE$ ( $= BF^2.CF$ ) de $BE^2.CF$:

$BE^2.CF - BE^2.CE = BE^2.CF - BF^2.CF$,  i.e

$BE^2. EF = (BE^2 - BF^2).CF$

et J. P. Hogendijk poursuit:

"*hence:*

$$BE^2 = CF.(BE + BF) \qquad (26)$$

*hence*

$$BF.(BE + BF - CB) = BE.CE. \qquad (27)$$

*In order to cast (27) in a nice geometrical form, al-Tûsî defines G on BC such that $BG = CE$ ( $= a - x_2$ ). Then (27) can be written as*

$$BF.GF = BE.CE. \qquad (28)$$

*If B, E and C (and hence G) are known, the construction of F is a standard Euclidien problem: to apply to BG a rectangle, equal in aerea to BE.CE, and exceeding*

*by a square (GF2). Or, in other words, GF2 + BG.GF = BE.CE (this is the equation used in [T2, 7] ). The fact that al-Tûsî uses GF and not BF (in (27)) as the unknown shows that his method is basically geometrical."* [5. p. 78. li. 16 - 26][6].

Ainsi, si al-Tûsî avait utilisé l'équation (27) (où l'inconnue est *BF*) au lieu de (28) (où l'inconnue est *GF*), on n'aurait pas pu qualifier son raisonnement de "géométrique". Or, la raison d'utilisation de (28) au lieu de (27) ne réside pas dans la volonté d'al-Tûsî d'utiliser certains théorèmes de la géométrie euclidienne; mais elle se dégage de l'étude de la première partie du "Traité". L'équation (27) est, en effet, une équation du second degré en *BF*

$$BF^2 - CE.BF = BE.CE$$

(27')  $y^2 - (a - x_2)y = x_2(a - x_2)$

du type (8):

(8)  $y^2 - by = c$

qui a une seule racine positive ($y = BF$)[7]. Mais, al-Tûsî la transforme en l'équation

(28)  $GF^2 + BG.GF = BE.CE$,
$x^2 + (a - x_2)x = x_2(a - x_2)$

qui est une équation du type (7)

(7)  $x^2 + bx = c$,

simplement par ce qu'il n'a jamais résolu directement les équations du type (8), mais indirectement, en les transformant en des équations du type (7), par le changement de variable : $x = y - b$ [2, t.1, éq. 8, p. CXLVI]. Or, en appliquant cette consigne qu'al-Tûsî a suivi déjà dans la première partie du traité, ($x = y - b$, i.e. $y = x + b$; $b = a - x_2$), l'équation (27'), plus haut, équivalente à (27), se transforme en l'équation:

$$(x + b)^2 - b(x + b) = x_2(a - x_2),$$

i.e

$$x^2 + bx = x_2(a - x_2),$$

i.e

$$x^2 + (a - x_2)x = x_2(a - x_2),$$

---

[6] Les numérotations (26), (27), (28) sont celle de J. P. Hogendijk, [5].
[7] Voir la note 1 plus haut sur la numérotation des équations.

qui n'est autre que l'équation (27) utilisée par al-Tûsî.

Cet exemple montre qu'une "reconstitution de l'analyse" doit toujours prendre en considération l'unité du projet d'al-Tûsî et par conséquent celle du "Traité". L'emploi de disciplines, de conceptions ou de techniques nouvelles dans la "deuxième partie" n'est pas dicté par une modification de l'objectif du "*Traité*" ou de son unité (voir la note complémentaire 1, vers la fin de l'article).

## 2. LA GEOMÉTRIE DANS LA DEUXIEME PARTIE DU "*TRAITÉ*".
### 2. 1. Caractère algébrique de la mathématique du "*Traité*".

Les remarques précédentes permettent de remettre en question le point d'appui d'une conjecture qui impute au raisonnement d'al-Tûsî le seul caractère géométrique. On doit cependant remarquer que l'invalidité de cet argument d'appui n'entraîne pas nécessairement celle de la conjecture elle-même; cette dernière pourrait, en effet, être vraie pour d'autres raisons, à découvrir.

Or, la lecture de "l'Introduction" de [2], assistée par celle du texte du "Traité", est capable de nous éclairer sur la géométrie du "*Traité*": sa nature, sa destination, la place qu'elle occupe dans sa mathématique et surtout son rôle dans l'algèbre d'al-Tûsî.

C'est, en effet, le mot "algèbre" qui semble être le terme le plus approprié pour désigner la mathématique du "*Traité*". Rappelons à ce sujet des faits admis par tous et qui reflètent le caractère algébrique de cette œuvre:

1) <u>Le titre</u>: "*Sur les équations*".

2) <u>Le thème</u>: Résolution des équations polynomiales de degré inférieur ou égal à 3.

3) <u>Le classement de ces équations</u>: non suivant le degré, (linéaires, planes, solides) mais selon d'autres critères algébriques et analytiques (classement, qui mériterait une étude à part).

4) <u>La terminologie</u> qui est celle de l'algèbre de l'époque : "*chose*" ou, autrement, "*racine*", (pour l'inconnue $x$), "*mâl*" (pour $x^2$), "*cube*" (pour $x^3$),..., enrichie par des termes nouveaux: "*nombre d'écart*", "*partie propre*",.... A ce propos, Rashed souligne bien l'effet négatif des handicaps de la terminologie de l'époque et du symbolisme en général, sur cette algèbre et sur son développement ultérieur.

5) Les techniques calculatoires:

a/ <u>Celles qui sont relatives aux lois de l'arithmétique élémentaire pour les nombres réels et pour les polynômes</u>. Remarquons qu'ici, les propriétés de l'addition et de la multiplication: associativité, distributivité de la multiplication par rapport à l'addition, simplification par un facteur additif ou multiplicatif non nul (nombre ou polynôme)... ne sont ni interprétées ni démontrées géométriquement comme chez Euclide.

b/ <u>Celles qui sont relatives aux résolutions d'équations</u>, notamment la technique de changement de variable affine ramenant une équation à une autre qui est déjà résolue. Cette technique maîtrisée admirablement par al-Tûsî, caractérise son algèbre. Elle joue un rôle important dans l'ordre adopté par al-Tûsî pour l'ensemble des 25 équations traitées. C'est la technique dominante dans la deuxième partie du "*Traité*", où elle se trouve appliquée deux fois à chacune de ses équations: une fois pour la détermination de la petite racine et une fois pour celle de la grande.

Ces évidences ne font qu'aiguiser notre curiosité pour comprendre les raisons qui ont permis d'avancer des affirmations telles que:

"*This second part (i.e. la 2° partie du "Traité") consists mainly of a sequence of very long proofs in Euclidian style.*" [5, p.70, li. -10 ],

ou aussi que :

" *The second part of the "Algebra" (i.e. la 2° partie sus-mentionnée) contains very little of what we could call algebra, i.e., direct manipulation of algebraic equations*" [5, p. 73, li. -4].

Ainsi, nous jugeons utile d'examiner la géométrie utilisée dans la deuxième partie du "*Traité*" pour essayer d'évaluer son rôle avec plus de précision.

**2. 2. Rôle de la géométrie dans la deuxième partie du "*Traité*".**

Les équations de troisième degré de la "première partie" du "*Traité*", admettent toujours une racine positive, au moins. On y trouve un chapitre géométrique dont le rôle consiste à:

a) démontrer l'existence d'une racine pour l'équation considérée;

b) représenter graphiquement cette racine à l'aide d'intersection de coniques[8].

Par contre, dans la, deuxième partie, un tel souci d'al-Tûsî est écarté, parce que:

1) l'existence de racines est traitée analytiquement ( cf.1.1.4) );

2) la détermination des racines d'une équation se fait en ramenant cette dernière (par un changement de variable affine, donc d'une façon purement algébrique) à une équation déjà résolue dans la première partie.

Les méthodes géométriques pour résoudre les équations algébriques ont survécu longtemps à l'époque médiévale. Mais, dans la deuxième partie du "Traité", il ne s'agit pas de méthodes géométriques. Les figures qui s'y trouvent ne risquent pas d'induire le lecteur en erreur (cf. note 7). Il est superflu de les analyser pour s'apercevoir qu'il s'agit plutôt de représentations graphiques que d'outils de raisonnement géométrique. Notons aussi qu'on ne rencontre, nulle part dans cette deuxième partie du "Traité", de techniques et d'outils de base propres au raisonnement géométrique euclidien: droite menée parallèlement à une autre, triangles superposables ou semblables, angles égaux, supplémentaires ou complémentaires .. Les exceptions présentées dans la note (8) sont les seules que nous avons pu relever dans toute la deuxième partie du "Traité"; on les rencontre au cours de l'étude de la seule équation (24), surtout lors de la localisation des racines de celle-ci. Elles confirment ainsi, le rôle marginal de la géométrie dans cette partie.

Les notes complémentaires (3) et (4) (cf. 5. plus loin) confirment la remarque déjà formulée par Rashed: "*les figures géométriques ne sont là que pour aider l'imagination*" [2. t. 1; p. XXVII]. On pourrait toujours, interpréter géométriquement les équations et les résultats algébriques; on pourrait développer et étendre certaines figures géométriques[9], mais on ne devrait jamais confondre le "raisonnement géométrique" et la représentation géométrique des expressions algébriques.

---

[8] Ce double rôle est clairement exposé par Rashed dans [2] où, de plus, l'accent est mis sur l'apport d'al-Tûsî dans ce domaine; dans [4], Houzel présente un aperçu détaillé et concis sur les innovations de ce dernier par rapport à Al Khayyâm.

[9] Par exemple, les fig 2. et 3. de [5] ne sont pas celles qu'al-Tûsî a utilisées lors de son étude de l'équation (25).

## 3. "RECONSTITUTION DE L'ANALYSE".

Les travaux effectués en vue d'une "reconstitution de l'analyse" dans le "*Traité*" montrent que cette tâche n'est pas simple, c'est-à-dire qu'elle est composée d'une série de "reconstitutions" propres à chacun des points du "*Traité*" où l'analyse est absente. Notons toutefois, que la part la plus pénible de cette tâche se trouve déjà accomplie par Rashed qui a reconstitué les tableaux relatifs aux algorithmes numériques et commenté la traduction de chaque paragraphe du texte. Après la publication des œuvres d'al-Tûsî, certaines études aboutissent à des résultats qui ne concordent pas avec l'ensemble des travaux établis dans ce domaine. Nous pouvons ainsi constater qu'une "reconstitution" de l'analyse devrait obéir à des exigences dont nous retenons quelques-unes (bien qu'elles semblent évidentes); nous croyons, en effet, que la négligence de ces exigences ou, tout au moins, le manque de rigueur dans leur application constitue la cause essentielle de ces divergences:

1/ La compréhension de l'ensemble du travail présenté par le "*Traité*" est nécessaire, même pour la reconstitution d'un détail aussi petit qu'il semble être. Nous insistons ici particulièrement sur l'unité de cette œuvre d'al-Tûsî (voir la note complémentaire 1). (Nous avons rencontré, (cf. 1.2) un exemple sur les risques que pourrait engendrer le traitement d'une des deux "parties" de cette œuvre indépendamment de l'autre).

2/ Bien que la démarche d'al-Tûsî soit uniforme (dans ses grandes lignes) dans sa résolution des équations (21),..,(25), il serait risqué de généraliser le résultat d'une de ces équations aux autres[10]. Il est surtout recommandé, dans toute "reconstitution" de l'analyse, de ne considérer aucune de ces équations comme modèle pour les autres[11].

---

[10] Par exemple, d'après [5, p. 80, li. -9 .. -7], l'origine des expressions de la forme *f(BD) - f(BE)* revient au paragraphe consacré à la démonstration de l'algorithme numérique propre à l'équation (21) "*and possibly for similar proofs for eqs. (2) - (5) - i.e. (22),..,(25) -*". Or, il n'y a pas d'algorithme numérique pour le calcul de la petite racine des équations (22), .., (25), vu que, pour le calcul de la petite racine, ces dernières équations sont ramenées à des équations du type (21).

[11] L'équation (25) présente des particularités spécifiques qui l'empêchent de constituer un modèle en vue d'une reconstitution de l'analyse [5]. En effet, al-Tûsî envisage une multitude de cas (par souci de ne rencontrer que des expressions algébriques positives, comme l'a déjà signalé R. Rashed). Cette équation occupe ainsi, une part importante dans la deuxième partie du "Traité"; al-Tûsî -probablement par souci de concision- a dû ne pas reprendre certaines notations et démarches qui se trouvent développées ailleurs, à des endroits homologues concernant les autres équations (cf. [6]). Disons enfin que l'itinéraire de cette reconstitution proposée dans [5], se démarque nettement de celui dessiné par al-Tûsî (par exemple,

3/ La découverte récente de l'œuvre d'al-Tûsî laisserait supposer, dans une lecture rapide, que telles ou telles interprétations sont optimistes ou excessives. Cependant, la présence du texte et de sa traduction d'une rigueur extrême, peuvent préserver le checheur des interprétations subjectives. Ainsi, le problème qui se pose à l'historien est de comprendre les méthodes et les motivations d'al-Tûsî pour les situer dans la diachronie mathématique, en se basant surtout sur le texte lui-même. Il ne consiste nullement à savoir "whether al-Tûsî's methods and motivation can also be explained in terms of standard ancient and medieval mathematics" [5. p. 71. li.3 - 4]. Il est possible, en effet, d'expliquer des mathématiques, même récentes, en "termes standard ...", sans qu'une telle explication apporte des indications sur leur positionnement historique ou sur leur classification mathématique. Notons surtout que l'ordre selon lequel les hypothèses et les résultats se présentent et se succèdent dans le "*Traité*" devrait être respecté.

Les deux détails que nous développons ci-dessous constituent des exemples sur les conséquences risquées de reconstituer une analyse sans prendre en considérations les consignes qu'on vient de citer. Ils contiennent surtout des indices significatifs sur le rôle du concept du "maximum" dans la démarche d'al-Tûsî.

**3. 1. La notion du maximum et l'existence des racines.**

Dans la démarche d'al-Tûsî pour résoudre les équations $f(x) = c$ de la "deuxième partie" du "Traité", l'existence et la détermination du maximum $f(x_0)$ précèdent l'existence et la détermination des racines, comme le montre le schéma 1.1. Al-Tûsî affirme sans aucune explication que, dans le cas où $c < x_0$, l'équation (E) en question, admet deux racines distinctes $x_1$ et $x_2$ telles que

$$x_1 < x_0 < x_2;$$

la détermination de $x_2$ (en passant à une équation du type (15) déjà résolue[12]) puis celle de $x_1$, suivront. C'est dans la phase finale qu'il calcule $x_1$ numériquement.

Nous n'interprétons pas, ici, les raisons qui ont mené al-Tûsî à énoncer l'existence des racines $x_1$ et $x_2$, à la suite de la détermination de $x_0$ (même si ces raisons semblent claires). Mais, nous attirons l'attention sur les risques que pourrait courir toute description

---

contrairement à l'ordre suivi par Al-Tûsî, le cas $a < b^{1/2}$ a été traité avant le cas $a > b^{1/2}$), ce qui ne peut se faire sans risquer des divergences du type rencontré lors de l'étude du maximum (cf. 3.1).

[12] La seule racine positive de l'équation (15) a été déterminée à l'aide de l'intersection d'une branche d'une hyperbole équilatère et d'une branche d'une parabole [2, t. 1, p. CLXIII].

de la démarche d'al-Tûsî, ne prenant pas en considération cet enchaînement précis. Considérons, par exemple, la description qui ne présente pas ce schéma comme la suite des pas: [1), 2)], 3), 4) (= 4.1, 4.2, 4.3), 5), 6), 7), 8) (cf. 1.1) mais comme la suite: [1), 2)], 3), 4) (= 4.1, 4.2), e'), ...où, dans e'), on suppose qu'al-Tûsî détermine $x_2$ en passant à une équation du type (15) pour en déduire, ensuite, son existence [5. p. 72. li. 24, 25].

Négliger la composante 4.3 de 4), mène à considérer l'existence de la plus grande racine, comme une conséquence de la résolution de l'équation (15) (celle-ci est établie géométriquement par al-Tûsî (5)). Or, dans 5), al-Tûsî ne prouve pas l'existence d'une racine $x_2 > x_0$, mais "détermine" cette racine; l'existence de $x_2$ semble être (pour lui) admise au préalable comme évidente, à la suite de la détermination du maximum $f(x_0)$. On ne trouve pas dans le texte se rapportant à l'étape 5), de phrase équivalente à "*donc $x_2$ existe*"; mais, en conclusion de sa détermination de $x_2$ il dit: "*Ce que nous cherchons dans ce problème est BE et il est plus grand que les 2/3 de AB* (i.e. il affirme seulement que $x_2 > 2a/3$) ". [2; t.2; p. 6, li.18]. Pourtant, al-Tûsî a été clair auparavant, en parlant de l'existence de $x_1$ et de $x_2$ (cf. 1.1.d) ): "*et, si c est plus petit -que $c_0 = f(x_0)$- l'équation admet deux solutions dont l'une est plus grande que le 2/3 des mâls -i.e. que 2a/3- l'autre est plus petite*" [2,t. 2, p. 5, li. 10].

Un tel éloignement, apparemment infime, de l'itinéraire suivi par al-Tûsî, ne conduit pas seulement à la confusion de la "détermination" et de "l'existence" d'une racine, mais il contribue aussi à masquer le rôle du "maximum" dans l'établissement de l'existence des racines et à attribuer ce rôle à des considérations algébrico-géométriques: changement de variable affine pour passer à une équation du type (15), suivi d'une détermination géométrique de la racine positive de celle-ci.

Il semble aussi qu'une lecture du schéma 1.1., suivant un ordre différent de celui du texte a conduit la même étude [5] à une autre divergence. Cette étude cherche, en effet, la découverte par al-Tûsî de la valeur $x_0$ (telle que $f(x_0)$ est maximum), là où cette découverte ne peut pas exister: dans un "passage obscur" consacré par al-Tûsî à la preuve de son algorithme numérique propre au calcul de la "petite racine" $x_1$ de l'équation (21) : "*Hence al-Tûsî may well have discovered the substitution  y = (2/3).a - x  ($z_2 = m - x_2$ in the notation of section 2) -i.e. $x_0 - x_1$, dans nos notations, qui sont celles de [2]- in connection whith his investigation  of the proof of the algorithme for eq.(1) -i.e. pour l'équation (21), selon le classement d'al-Tûsî-" [5, p. 81, li. 5-7]. Elle laisse donc supposer que la découverte de $x_0$ est occasionnelle et qu'al-Tûsî ne l'aurait pas découvert en cherchant le maximum intentionnellement. Cette même étude attribue au dernier passage, qualifié à juste titre d'"obscur", l'origine des expressions algébriques de la forme $f(v) - f(u)$ ($f(BD) - f(BE)$), où $BD = x_0$ ( $= m$ ), $BE = x_2$, $BF = x_1$, cf. fig. 1): "*Hence it is conceivable that al-Tûsî first studied the differences f(BD) - f(BE) while he was searching for this proof, and possibly for similar proofs for eqs. (2) - (5) - i.e. (22)..(25) -*" [5; p. 80. li. - 9]. Or, ces expressions $f(v) - f(u)$ "*play a cardinal role in the reasoning of al-Tûsî*" [5; p. 78 - 79, li. 1], d'où le rôle capital de ce passage obscur dans la mathématique d'al-Tûsî. Nous jugeons ainsi utile de présenter une lecture concise de ce passage (voir 3.2, plus bas). Mais, précisons d'abord, que la détermination de $x_0$, le calcul du maximum, la détermination et le calcul de la plus grande racine ($x_2 = x_0 + y$), la détermination de $x_1$ lui-même ($x_1 = x_0 - y$) constituent des paragraphes qui ont utilisé $x_0$ comme concept et comme valeur effective. Ces paragraphes ont bien précédé, dans la démarche d'al-Tûsî, le passage en question. De plus, dans tous ces paragraphes apparaissent des comparaisons entre $f(v)$ et $f(u)$ et, par suite, des expressions de la forme $f(v) - f(u)$ (dans le cas où une telle différence est positive) et cela est naturel d'après le schéma 1.1.

**3. 2. Lecture du passage concernant le calcul de $x_1$.**

Nous devons avouer que la lecture de ce passage aurait été une tâche très difficile sans la traduction et les commentaires présentés dans [2], vu le langage mathématique qui y est utilisé par al-Tûsî. Un des objectifs visés par celui-ci, est de démontrer qu'il peut utiliser un algorithme unique pour le calcul numérique de la petite racine, $x_1$. Cet objectif risque de rester dans l'ombre si on lit ce passage indépendamment d'un lemme qui le précède de quelques pages (cf. lemme 2 [2, t. 2, p. XVII ou p. 11]). Dans ce lemme qui mériterait une étude à part, Al-Tûsi remarque que dans l'équation (21): $ax^2 - x^3 = c$, on a trois cas à considérer suivant les valeurs de c:

1° cas: $c = c_0/2 = f(x_0)/2$, alors $x_1 = a/3$,
2° cas: $c > c_0/2 = f(x_0)/2$, alors $x_1 > a/3$,
3° cas: $c < c_0/2 = f(x_0)/2$, alors $x_1 < a/3$.

(Signalons le rôle essentiel de $x_0$ ($= m$) dans ce lemme qui a précédé, de loin, l'algorithme calculatoire de $x_1$ [1. p.15 - 18]).

De ce "lemme 2", al-Tûsî déduit que:

1/ Si $c = c_0/2$, alors le calcul de $x_1$ est évident;
2/ Si $c < c_0/2$, alors $a - 3x_1 > 0$ et il peut appliquer son algorithme;
3/ Si $c > c_0/2$, alors $a/3 < x_1 < 2a/3$, alors il considère $y = 2a/3 - x_1 < a/3$; dans ce cas, $y$ est racine d'une équation du type (21) à laquelle il applique le même algorithme (vu que $y < a/3$ ). Le changement de variable $y = 2a/3 - x$ est, dans cet algorithme, dicté par la volonté d'al-Tûsî de démontrer que son algorithme s'applique bien aux différents cas qui se présentent ($c < c_0/2$, $c = c_0/2$ et $c > c_0/2$); la découverte de $x_0 = 2a/3$ étant préalable à cet algorithme.

### 3. 3. Al-Tûsî et le calcul du maximum.

Le calcul du maximum de $f(x)$ est essentiel dans toutes les étapes de la démarche d'al-Tûsî pour résoudre les équations (21), .. ,(25) (cf. 1.1.). Pourtant, il avance la valeur effective $v$ de $x_0$ (valeur de x pour laquelle f(x) est maximum), puis il démontre que $f(v)$ est maximum, sans toutefois indiquer explicitement les chemins par lesquels il est arrivé à trouver $v$. Son style concernant ce point est essentiellement synthétique (cf. 0.,1)). Une reconstitution de son analyse s'impose donc; elle consiste à découvrir les moyens qui l'ont conduit à trouver cette valeur $v$ et par suite à annuler l'expression de la dérivée de certaines expressions polynomiales. Car, bien que, pour les équations (21) et (22) al-Tûsî donne $v$ par sa valeur ($v = 2a/3$ et $v = \sqrt{b/3}$, resp.), pour les équations (23), (24) et (25), il présente $v$ explicitement comme une "racine d'une équation" équivalente à $f'(x) = 0$ (plus précisément $f'(x)/3 = 0$; probablement le facteur multiplicatif 1/3 a été introduit pour présenter cette équation sous sa forme normalisée, conformément à la tradition de l'époque).

Cette reconstitution de l'analyse est donc capitale. Sa réalisation d'une façon définitive aiderait à mieux caractériser cette mathématique du XXII[ème] siècle et à évaluer sa portée, surtout dans le domaine de l'analyse mathématique. A ce propos, R. Rashed présente une conjecture basée sur une étude du texte d'al-Tûsî, (cf. R. Rashed [2. t. 1, p. XXIV]; cf. aussi, C. Houzel [3]). D'autre part, nous venons dans le paragraphe précédent, de discuter une interprétation concernant le même sujet. Or celle-ci, comme on vient de le voir, s'appuie sur l'équation (21), la seule dans la deuxième partie du "*Traité*" pour laquelle al-Tûsî donne un algorithme calculatoire de la petite racine. Cette interprétation ne peut donc, même si on l'admet, proposer une réponse générale concernant les quatre autres équations du troisième degré (22),..,(25) de la deuxième partie du "Traité".

Enfin, dans une tentative de contribution à la reconstitution de l'analyse relative à ce point, nous croyons avoir démontré qu'al-Tûsî est arrivé à la valeur *v* de $x_0$ au cours de sa démonstration des relations

(1) $$\begin{cases} x < x_0 \Rightarrow f(x) < f(x_0) \\ x > x_0 \Rightarrow f(x) < f(x_0), \end{cases}$$

qui constituent pratiquement une définition du maximum de *f*. Cela veut dire qu'il a calculé intentionnellement *v* en partant de la définition (1) de $x_0$; laquelle définition se trouve rappelée par al-Tûsî au début de l'étude de chacune des équations (21),..,(25). Notre travail concernant ce sujet [6], bien qu'il soit différent (surtout par son point de départ) de celui de Rashed, ne contredit pas la conjecture de celui-ci; Il s'harmonise plutôt avec sa conclusion affirmant que la présence de certaines notions qui appartiennent au domaine de l'analyse mathématique "*n'est ni fortuite ni secondaire mais plutôt intentionnelle*". Notons que l'hypothèse selon laquelle al-Tûsî n'a pas d'abord calculé la "fonction dérivée": "*al-Tûsî does not find m -i.e. $x_0$- by computing the derivative f' and by puting f'(x) = 0*" [5. p. 75. li.-5], semble être incontestable; on était bien loin de l'époque où la notion de "fonction" s'était cristallisée. Mais cela n'empêche qu'il est arrivé à annuler l'expression de la dérivée en un point. La définition (1) du maximum conduit naturellement à l'équation *f'($x_0$) = 0*, pour toute fonction *f* dérivable en $x_0$.

**4. CONCLUSION.**

Notre question de départ portait sur la possibilité de classer la mathématique du "*Traité*" et surtout celle de sa "deuxième partie" dans le seul domaine de la géométrie euclidienne. Pour y répondre nous avons rappelé que le travail d'al-Tûsî s'inscrit dans une stratégie essentiellement algébrique. Nous avons surtout insisté sur:

1°/ le rôle de certaines techniques algébriques (notamment le changement de variable affine) dans la résolution des équations et dans la présentation de leur nomenclature,

2°/ le rôle des notions analytiques comme celle du maximum d'une expression polynomiale $f(x)$ qui, curieusement, a conduit ce mathématicien du douzième siècle à des équations équivalentes à $f'(x)=0$.

Cette mise au point montre que la deuxième partie du "*Traité*" ne peut pas être comprise profondément à partir de la seule géométrie euclidienne: il s'agit, en fait, d'une étude essentiellement algébrique, comprenant également des méthodes analytiques.

**5. NOTES COMPLEMENTAIRES**.
**(1)** Le passage suivant extrait de [2, t. 1, pp. XXXI-XXXII], rappelle d'une façon claire et concise l'unité du projet -donc celle du "Traité"- d'al-Tûsî. Il donne une idée de la différenciation d'ordre mathématique de ces deux parties, le rôle de la géométrie dans la première et l'émergence des conceptions analytiques dans la deuxième*:*

*"Reconstituer l'unité d'un chapitre réservé à la théorie des équations algébriques, telle est, selon nous, la visée qui domine le "Traité"; et la rédaction d'al-Tûsî demeure obscure tant que n'est pas suffisamment dégagée cette intention délibérée d'élaborer un exposé structuré et unitaire. Et pourtant, c'est précisément ce projet qui n'a pu résister à la construction du "Traité": l'unité voulue est en effet brisée par l'émergence, au cours de la recherche, d'une problématique que rien ne laissait prévoir initialement, et qui scinde le "Traité" en deux parties, solidaires certes, mais relevant de deux mathématiques différentes. La première partie, dans la tradition d'Al Khayyâm, est destinée à la construction géométrique des racines des équations; c'est au cours de cette recherche qu'al-Tûsî s'impose une contrainte supplémentaire: démontrer systématiquement, ce qui veut dire pour chaque cas, l'existence du point d'intersection de deux courbes dont*

*l'abscisse détermine la racine positive demandée. Cette nouvelle exigence a, tout naturellement, conduit l'auteur à poser les problèmes de localisation et de séparation des racines, et à traiter des condition de leur existence indépendamment de leur constructibilité géométrique. C'est pour résoudre ce problème qu'al-Tûsî définit la notion de maximum d'une expression algébrique, et s'efforce de trouver concepts et méthodes pour la détermination des maxima. Non seulement le mathématicien est-il amené, dans cette démarche à inventer notions et méthodes qui ne seront baptisées que plus tard, mais il doit, pour y parvenir, changer de mode d'approche: pour la première fois, que je sache, il découvre la nécessité de procéder localement. La seconde partie du "Traité", précisément consacrée à ces problèmes, se différencie de la première par les objets qu'elle considère, et s'en démarque plus par les style mathématique qu'elle adopte."*

**(2)** Les 79 figures géométriques de la deuxième partie du "*Traité*" peuvent, en effet, être classées dans l'un ou l'autre des cinq types de figures:

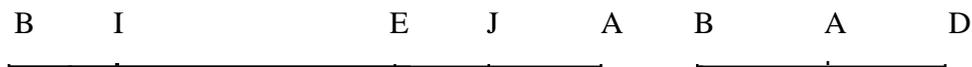

(fig. type 1)

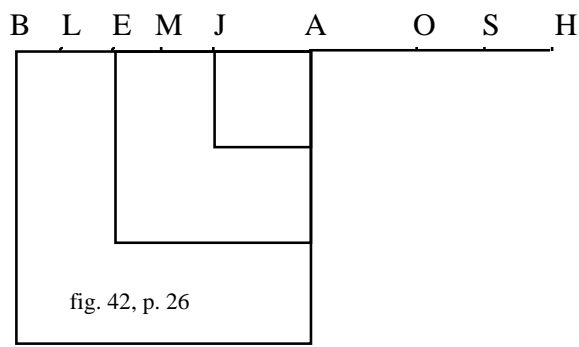

fig. 42, p. 26

(fig. Type 2)

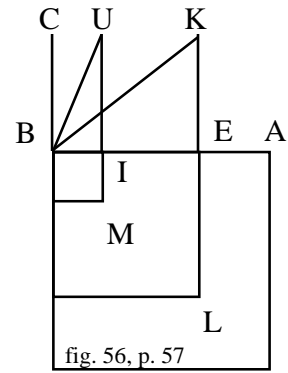

fig. 56, p. 57

(fig. Type 3)

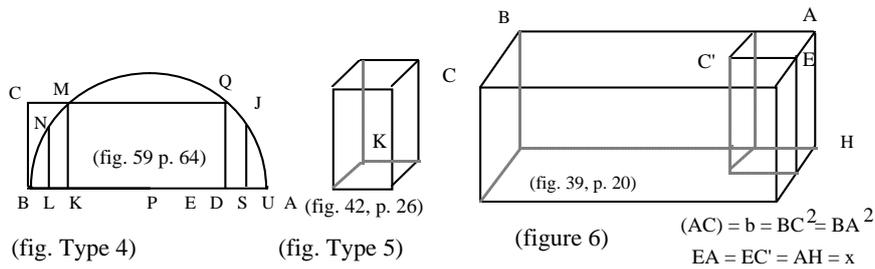

(fig. Type 4)  (fig. Type 5)  (figure 6)

$(AC) = b = BC^2 = BA^2$
$EA = EC' = AH = x$

Elles sont réparties de la façon suivante:

- 51 figures du type 1, réparties sur les cinq équations (21),..,(25), où des racines et des coefficients de l'équation étudiée, sont représentés sur un même segment de droite approprié.

- 15 figures du type 2 (dont on présente ici l'une des formes les moins simples) où BA est du type 1. Les figures de ce type, fréquentes surtout dans les équations (22) et (23), sont présentées souvent sans le segment AH ou sans l'un ou l'autre des carrés intérieurs.

- 6 figures du type 3 que l'on rencontre dans l'équation (24) lors de la détermination de $x_0$, $x_1$ et $x_2$. (BA étant du type 1). On présente ici (cf. plus haut) la forme la moins simple des figures de ce type ( cf [2. t. 2, p. 57, fig 56] ).

- 3 figures utilisant des demi-cercles, données à l'occasion de l'encadrement des racines de l'équation (24) (cf. [2. t. 2, fig.53, 59, 60; pp. 50, 64, 66]). Nous pouvons les grouper (à une nomination de points près) en une seule figure (fig. type 4).

- 3 figures où l'on rencontre des parallélépipèdes sans aucune construction géométrique sur aucun d'eux.( fig type 5 ).

- Une unique figure dans l'espace représentant $bx - x^3$, sans aucune autre construction géométrique là-dessus (fig. 6), où $BC^2 = BA^2 = (AC) = b$;  $EC' = EA = AH = x$.

**(3)** Pour plus de précision, notons les trois exceptions suivantes, rencontrées lors de l'étude de la seule équation (24):

1°) Deux segments de droite *IU* et *EK* ont été menés perpendiculairement à AB (fig. type 3). Un demi-cercle de centre *P* et de diamètre *AB* a été tracé, ainsi que deux segments de droite *MK* et *JU* perpendiculaires à *AB* (fig. type 4).

2°) Le théorème de Pythagore a été appliqué deux fois (cf. fig. type 3) aux triangles *BEK* et *BIU* ([2. t.2, p. LXII et p. LXIII; p. 56, li. 5-6 et p. 57, li. 7 ), puis rappelé, une fois, ultérieurement (p. LXXI; p. 68, li. 15).

3°) La puissance d'un point par rapport à un cercle ($KM^2=BK.KA$) a été utilisée (cf. fig. type 4, - p. LIX, p. 50, li. 4 -), lors de la démonstration de l'impossibilité de résoudre l'équation (24) si $b^{1/2} > a/2$, de même que lors de l'encadrement des racines (p. LXVII; p. 63, li. 9). On trouve aussi utilisée une relation qui revient à $BU.AU < SQ^2$ si *U* est un point du segment *SA* (p. LXIX; p. 65. li 18-19).

******************************************

## BIBLIOGRAPHIE

[1] R. Rashed: "Résolution des équations numériques et algèbre. Sharaf al-Dîn al-Tûsî, Viète"; *Archive for History of Exact Sciences*; vol. 12; n° 3; 1974, pp. 244 - 290; cet article a été reproduit dans le chap. 3, pp. 148 - 193, du livre du même auteur: "*Entre arithmétique et algèbre: Recherches sur l'Histoire des mathématiques arabes*", Paris, Belles Lettres, 1984.

[2] R. Rashed: "*Sharaf al-Dîn al-Tûsî: Œuvres mathématiques. Algèbre et géométrie au XIIe siècle*". T. 1 et 2. Paris, "Les belles lettres", 1986. Traduit en arabe par N. Farès: éd. "Centre of Arab Unity Studies" Beyrouth (A paraître) 1995.

[3] C. Houzel: "Œuvres matématiques: Algèbre et géométrie au XIIème siècle; Sharaf al-Dîn al-Tûsî": Compte-rendu du livre du même titre; *Gazette des mathématiciens n° 39*. Paris, Janvier 1989, pp.59 - 63.

[4] C. Houzel: "Algèbre et géométrie au XIIème siècle; Sharaf al-Dîn al-Tûsî". *Arabic Sciences and philosophy. Cambridge University Press*. Vol. 5.2. 1995.

[5] J. P. Hogendijk: "Sharaf al-Dîn al-Tûsî: On the number of positive roots of cubic equations". *Historia mathématica*, 16 (1989), pp. 69-85.

[6] N. Farès: "Le calcul du maximum et la "dérivée" selon sharaf al-Dîn al-Tûsî". *Arabic Sciences and philosophy. Cambridge University Press. Vol. 5.2. 1995*.

[7] R. Rashed et B. Vahabzadeh "*Al-Khayyâm mathématicien*", Blanchard, Paris, 1999.